\documentclass[a4paper,oneside,11pt]{article}
\usepackage[english]{babel}
\usepackage[latin1]{inputenc}
\usepackage[T1]{fontenc}
\usepackage{amsmath}
\usepackage{amsfonts}
\usepackage{amssymb}
\usepackage{dsfont}
\author{Gabriel Vigny}
\begin{document}
\newtheorem{theorem}{Theorem}[section]
\newtheorem{proposition}[theorem]{Proposition}
\newtheorem{defi}[theorem]{Definition}
\newtheorem{corollaire}[theorem]{Corollary}
\newtheorem{lemme}[theorem]{Lemma}	
\title{Dynamics semi-conjugated to a subshift for some polynomial mappings in $\mathbb{C}^2$}
\maketitle	

\begin{abstract}
 We study the dynamics near infinity of polynomial mappings $f$ in $\mathbb{C}^2$. We assume that $f$ has indeterminacy points and is non constant on the line at infinity $L_\infty$. If $L_\infty$ is $f$-attracting, we decompose the Green current along itineraries defined by the indeterminacy points and their preimages. The symbolic dynamics that arises is a subshift on an infinite alphabet. 
\end{abstract}

\noindent\textbf{MSC:} 37F10, 37B10\\
\noindent\textbf{Keywords:} Polynomial mapping, Green current, subshift.

\section{Introduction}
  
  We are interested in the dynamics of polynomial mappings $f$ in $\mathbb{C}^2$ whose meromorphic extensions to $\mathbb{P}^2$ admit indeterminacy points and for which the line at infinity (which we denote by $L_{\infty}$) is $f$-attracting (that is: there exists $C>1$  such that for $p\in\mathbb{C}^2$ with $\| p \|$ large enough, one has $\| f(p) \| \geq C \| p \|$). In particular, given any large ball $\mathbb{B}$ in $\mathbb{C}^2$, these maps are polynomial-like in the sense of \cite{DS3} from $f^{-1}(\mathbb{B})$ to $\mathbb{B}$. The dynamics is studied there: there exists an invariant probability measure which is K-mixing and of maximal entropy. Our goal is to study the dynamics near infinity, especially the structure of the \emph{Green current}, which is a positive closed current of bidegree (1,1) invariant under the action of $f^*$. \\
  
In \cite{DDS}, the authors consider the case where $f_{\infty}$, the restriction of $f$ to $L_{\infty}$, is constant and they decompose the Green current into pieces associated to an itinerary defined by indeterminacy points. On the basin of attraction of the indeterminacy set, the itinerary map semi-conjugates $f$ to a shift. \\
    
    Another case which has been studied is when $f$ admits a holomorphic extension to $\mathbb{P}^2$: in \cite{BJ1}, the authors showed the Green current admits a local laminar decomposition consisting of local stable manifolds of $f$ to the Julia set of $f_{\infty}$. Applying one dimensional theory, one also obtains in this case a dynamics semi-conjugated to a shift. \\
    
    We study here a mixed situation. We assume that $f$ admits indeterminacy points on $L_\infty$ and that $f_\infty$ is not a constant function. In order to describe clearly the new phenomena happening here, we consider the case where $f_\infty$ is hyperbolic. The method we use allows to study more general cases. We will complete our study by giving several examples. In the hyperbolic case, we show that the Green current decomposes along some itineraries defined by the indeterminacy points and theirs preimages. Surprisingly, the local stable manifolds associated to the Julia set of $f_\infty$ are not charged by the Green current. Furthermore, the symbolic dynamics we obtain is a subshift (a Markov chain), which is new for polynomial mappings.\\
    
   The main tools we use are horizontal-like maps  and a theorem of convergence of currents proved in \cite{DUJ} and \cite{DDS}. Roughly speaking, such applications are contracting in the vertical direction and expanding in the horizontal one in a geometrical sense. For the reader's convenience, we give the basic properties of those objects. 
 
 Next, we define and study the basic properties of the family $\mathcal{G}$ of maps we consider. We give a simple sufficient condition for a map $f$ to be in $\mathcal{G}$ and we prove the algebraic stability. Then, by a theorem of Sibony \cite{SIB}, one can associate to $f$ a natural invariant current (Green current). We give an easily computable formula for the trace of the Green current at infinity. This trace is a probability measure which is a combination of Dirac masses at the indeterminacy points and theirs preimages. Under some additional hypothesis, we also compute the topological degree.
    
    We then study the decomposition of the Green current on a neighborhood of infinity under the hypotheses that the indeterminacy set is located in the Fatou set of $f_\infty$, with no indeterminacy point being periodic for $f_\infty$ and that $f_\infty$ is hyperbolic. This set of maps contains an open subset of $\mathcal{G}$. The decomposition of the Green current semi-conjugates $f$ to a subshift on an infinite alphabet. Under some additional hypothesis, we show that the range of the \emph{escape rate} (which measures the asymptotic speed at which a point goes to infinity) is a full interval which is new for polynomial maps and we compute a mean escape rate. We will explain briefly how to obtain a weaker decomposition of the Green current in a more general case. Finally, we study examples, in particular the case where the indeterminacy points are located in the exceptionnal set of $f_{\infty}$, in this case the support of the Green current is strictly contained in the Julia set of $f$. \\

\section{Polynomial maps with dynamics at infinity}
\subsection{Horizontal-like maps}\label{tool} 
 We recall here the facts we use on horizontal-like maps. Proofs and details can be found in \cite{DUJ} and \cite{DDS}. 
 
 Let $\mathbb{D}$ (resp. $\mathbb{D}_r$) be the unit disk  (resp. the disk of radius $r$ centered at $0$) in $\mathbb{C}$. Let $\Delta$ be the unit bidisk in $\mathbb{C}^2$, we denote its vertical boundary by $\partial_v \Delta$, and its horizontal boundary by $\partial_h \Delta$. Namely:
 $$\partial_v \Delta  =   \{(z,w)\in \mathbb{C}^2, \ |z|=1,\ |w|<1 \} \ \text{and} \ \partial_h \Delta =   \{(z,w)\in \mathbb{C}^2, \ |z|<1,\ |w|=1 \}.$$ 
 We have the following definitions: 
\begin{defi} Let $\Delta_i \subset M_i$ be an open subset biholomorphic to $\Delta$ in the complex surface $M_i$ for $i=1,2$. Let $f$ be a dominating meromorphic map defined in some neighborhood of $\Delta_1$ with values in $M_2$. The triple $(f,\Delta_1,\Delta_2)$ defines a horizontal-like map if:
\begin{itemize}
\item $f$ has no indeterminacy points in $\partial_v \Delta_1$ and $f(\partial_v \Delta_1) \cap \overline{\Delta_2}=\emptyset$;
\item $f(\overline{\Delta_1}) \cap \partial \Delta_2 \subset \partial_v \Delta_2$;
\item $f(\Delta_1)\cap \Delta_2 \neq \emptyset$.
\end{itemize}
 \end{defi}
\begin{defi}A positive closed (1,1)-current $T$ in $\Delta$ is \emph{vertical} if:
$$ \mathrm{Supp T} \subset \mathbb{D}_{1-\varepsilon}\times \mathbb{D} \ \text{for some} \ \varepsilon>0.$$
\end{defi}
Similarly, we can define horizontal currents.\\

We can define the (horizontal) slice measures $m^{w_0}$ of a vertical positive closed (1,1)-current $T$ by $T \wedge [w=w_0]$. These measures have the same mass, which we call the \emph{slice mass} of $T$. The current $T$ is zero if and only if its slice mass is zero. The main fact is that we can define the pull-back of such a current by a horizontal-like map, and we have the following:\\

\emph{Let $(f,\Delta_1,\Delta_2)$ be a horizontal-like map then there exists a positive integer $d\geq 1$ such that for every vertical positive closed current $T$ in $\Delta_2$ of slice mass $1$, $\frac{1}{d} f^{*}(T)$ is a vertical positive closed current in $\Delta_1$ of slice mass $1$.}\\

We call this integer the degree of $f$, it can be computed as the intersection multiplicity of the preimage of a vertical line with a horizontal line. The following result is our main tool to obtain the convergence in Theorem \ref{main}:
\begin{theorem}[\cite{DDS}]\label{convergence} Let $\{(f_i,\Delta_i,\Delta_{i+1} )\} _{i\geq 1}$ be a sequence of horizontal-like maps of degree $d_i$ such that $(f_i)^{-1}(\Delta_{i+1})\subset \mathbb{D}_{1-\varepsilon} \times \mathbb{D} \subset \Delta_{i}$ for a fixed $\varepsilon>0$. 
Assume that $K=\bigcap_{n\geq 1} f_1^{-1}\dots f_n^{-1}(\Delta_{n+1})$ has zero Lebesgues measure. For each $n$, let $T_n$ be a vertical positive closed (1,1)-current of slice mass $1$ in $\Delta_n$.

Then, the sequence of iterated pull-back $(\frac{1}{d_1}f_1^{*} \dots \frac{1}{d_n}f_n^{*} T_{n+1})_n$ converges to a vertical positive closed current $\tau$ of slice mass $1$ in $\Delta_1$ which is independent of $(T_n)$.  
\end{theorem}

\subsection{The class $\mathcal{G}$}
 We are interested in the study of the polynomial mappings $f$ of $\mathbb{C}^2$ for which $L_{\infty}$ is attractive. Namely, there are constants $C >1$ and $M$ large enough such that for $\|p\| \geq M$, we have $\|f(p)\| \geq C\|p\|$. We also assume that the meromorphic extension of $f$ to $\mathbb{P}^2$ admits indeterminacy points, we still denote by $f$ that extension. The case where the restriction $f_{\infty}$ of $f$ to $L_{\infty}$ is constant was studied in \cite{DDS}, so we will assume that $f_\infty$ is not constant. We denote by $\mathcal{G}$ the set of mappings satisfying the above properties (in fact, we are particularly interested in the case where $f_\infty$ is hyperbolic and where the indeterminacy points are not periodic for $f_\infty$). \\
 
 Let $f=(f_1,f_2)$ of algebraic degree $D$ be in $\mathcal{G}$. We denote by $f^+_1$ and $f^+_2$ the homogeneous parts of maximal degree. After a linear change of coordinates, we can assume $\deg f^+_1=D$ and $\deg f^+_2=D'\leq D$. The meromorphic extension of $f$ to $\mathbb{P}^2$ is given by $[T^D f_1(Z/T,W/T):T^D f_2(Z/T,W/T):T^D]$ and the restriction of $f$ to $L_{\infty}=(T=0)$ is $f_{\infty}[Z:W]=[f^+_1(Z,W):0^{D-D'}f^+_2(Z,W)]$. Thus, in order to have $f_{\infty}$ not constant, we need $D=D'$ and $f^+_1$ not proportional to $f^+_2$ (otherwise, $f$ sends $L_{\infty}$ to $[1:0:0]$ or $[1:\lambda:0]$).\\
 
 The indeterminacy set $I(f)$ of $f$ is the common zeros of $f^+_1$ and $f^+_2$: if the line $\mathcal{D}$ of equation $a_j z - b_j w =0$ satisfies $f^+_1(\mathcal{D})=\{0\}$ and $f^+_2(\mathcal{D})=\{0\}$ then $[b_j:a_j:0]$ is in $I(f)$.\\
 
 One deduces from above that all the mappings of $\mathcal{G}$ can be written as:
 \begin{eqnarray}
 \label{form}
 f(z,w) &=& \Big(\prod_{j\leq m} (a_j z- b_j w)^{\alpha_j} P_1(z,w) + Q_1(z,w),\nonumber \\
            & & \ \prod_{j\leq m}  (a_j z- b_j w)^{\alpha_j} P_2(z,w) + Q_2(z,w) \Big), 
 \end{eqnarray} 
where the $a_j$ and $b_j$ are complex numbers satisfying $(a_j, b_j) \neq (0,0)$, $m$ and the $\alpha_j$ are positive integers, $P_1$ and $P_2$ are homogeneous polynomials of degree $d'\geq 1$ with no common factor and the $Q_j$ are polynomials of degree strictly smaller than the degree of $f$. We denote by $d$ the sum $\sum_{j\leq m} \alpha_j$, so that $f$ is of degree $d+d'=D$.\\

We have the following formula for the extension of $f$ to $\mathbb{P}^2$: 
 \begin{eqnarray*}
 f([Z:W:T]) &=& \Big[ \prod_{j\leq m} (a_j Z- b_j W)^{\alpha_j} P_1(Z,W)+T^{D}Q_1(\frac{Z}{T},\frac{W}{T}) :\\
               & & \ \prod_{j\leq m} (a_j Z- b_j W)^{\alpha_j} P_2(Z,W)+T^{D}Q_2(\frac{Z}{T},\frac{W}{T}): T^{D} \ \Big]
\end{eqnarray*}  
thus $f_{\infty}([Z:W])=[P_1(Z,W):P_2(Z,W)]$. Recall that the multiplicity of an indeterminacy point $I$ is the intersection multiplicity at $I$ of $L_\infty$ and $f^{-1}(L)$ where $L$ is a generic line. The indeterminacy points of $f$ are the $I_j=[b_j:a_j:0]$ with multiplicity $\alpha_j$. We assume of course that the $(a_j,b_j)$ are not proportional.\\

The following proposition shows that we can find $f$ with any given set of indeterminacy points with multiplicity and any given restriction at infinity. Furthermore, it shows that for $D \geq 3$, $\mathcal{G}$ correponds to a Zariski open set of the space of parameters of (\ref{form}). If $d=d'=1$, we may have to multiply $f$ satisfying the criterion below by a large enough constant in order to have that the infinity is attracting. 
\begin{proposition}\label{criterion}
 Let $f=(f_1,f_2)$ be as in (\ref{form}). Assume that the polynomial $\Phi= f_1 P_2-f_2 P_1$ is of degree $\geq 2+d'$. If for all $j$, $a_jz-b_jw$ does not divide the homogeneous part of maximal degree of $\Phi$, then $L_{\infty}$ is $f$-attracting. 
\end{proposition}
\emph{Proof.}
Let $f$ be as above  and $N$ be a small neighborhood  of infinity. Observe that for any neighboorhood $V$ of $I(f)$, there exists a constant $C$ such that if $p=(z,w) \in N\backslash V$ we have $C \|p\|^{D}\leq \|f(p)\|$. So we just have to prove the estimate on $V$. Since $P_1$ and $P_2$ have no common factor, there is  $\lambda>0$ such that $\max(|P_1(z,w)|,|P_2(z,w)|) \leq \lambda \|(z,w)\|^{d'}$ on $N$. The hypothesis implies that $|\Phi(z,w)|\gtrsim \|(z,w)\|^{\deg \Phi}$ near $I(f)$, hence:
$$
 2 \|f(z,w) \| \geq \frac{|\Phi(z,w)|}{ \max(|P_1(z,w)|,|P_2(z,w)|)}\gtrsim \|(z,w)\|^{2}. 
  $$
The proposition follows. $ \Box$ \\ 
 
Observe that for a generic map $g\in \mathcal{G}$, we have $\deg \Phi =2d'+d-1$. The criterion is not optimal, but it is generic for $D \geq 3$ and easy to check. If $D=2$, we obtain in the same way that $\|f(z,w) \| \gtrsim \|(z,w)\|$. 

We use the notation of (\ref{form}) in the following proposition.
 
 \begin{proposition}\label{composition}
 Let $f$ and $g$ be in $\mathcal{G}$ then $f \circ g \in \mathcal{G}$. More precisely, if $f=(PQ_1+R_1,PQ_2+R_2)$ is of degree $D$ and $g=(P'Q'_1+R'_1,P'Q'_2+R'_2)$ is of degree $D'$ then $f \circ g = ( P''Q''_1+R''_1,P''Q''_2+R''_2)$  where:
 \begin{eqnarray*}
 P'' = (P')^D P(Q'_1,Q'_2),\quad Q''_1  =  Q_1(Q'_1,Q'_2),\quad \text{and}\quad  \  Q''_2  =  Q_2(Q'_1,Q'_2).
 \end{eqnarray*}  
In particular, $f \circ g$ is of degree $D+D'$.
 \end{proposition}
  \emph{Proof.} With the above notations, the homogeneous part of maximal degree of the components of $f \circ g$ are equal to:
  $$ P(P'Q'_1,P'Q'_2)Q_1(P'Q'_1,P'Q'_2) = (P')^D P(Q'_1,Q'_2) Q_1(Q'_1,Q'_2)$$
and
 $$  P(P'Q'_1,P'Q'_2)Q_2(P'Q'_1,P'Q'_2) = (P')^D P(Q'_1,Q'_2) Q_2(Q'_1,Q'_2).$$
We only have to check that $Q_1(Q'_1,Q'_2)$ and $Q_2(Q'_1,Q'_2)$ have no common factor: else, since two homogeneous polynomials have no common factor if and only if they have no common non trivial zero and since $Q_1$ and $Q_2$ have no common factor, we would have that $Q'_1$ and $Q'_2$ have a non trivial common zero. $\Box$ \\
 
Recall that a meromorphic mapping $f:\mathbb{P}^2 \to \mathbb{P}^2$ is said to be \emph{algebraically stable} if no algebraic curve is sent to an indeterminacy point after some iterations, equivalently, if $f$ is of algebraic degree $D$ then $f^n$ is of degree $D^n$ for all $n\geq1$. For such maps, the \emph{Green function} and the \emph{Green current} are defined  by $G(p)=\lim_{n\to\infty}(\frac{1}{D^n} \log^+ \| f^n(p) \|)$ and $T=dd^cG$ (see \cite{SIB}). Moreover, if $S$ is a smooth positive closed (1,1)-current of mass 1 on $\mathbb{P}^2$ then $\frac{1}{D^n}(f^{n})^{*}(S) \to T$ in the sense of currents. We deduce the important following corollary from the previous proposition:
\begin{corollaire}\label{stability}
 Let $f\in \mathcal{G}$ and $n\in \mathbb{N}$, then $\deg f^n=D^n$ so $f$ is algebraically stable. Furthermore, $(f_{\infty})^n=(f^n)_{\infty}$.   
\end{corollaire}

\subsection{Multiplicity of the indeterminacy points, trace of the Green current at infinity} 
  Let $E$ denote the set $\bigcup_{n\geq0} f^{-n}(I(f))=\bigcup_{n\geq0} I(f^n)$. For $p\in E$, we denote by $\lambda_{p,n}$ the real number equal to the multiplicity at $p$ of $f^n$ as an indeterminacy point divided by  $D^n$, that is: $\lambda_{p,n}=\frac{\text{mult}_p(f^n)}{D^n}$ (these numbers will appear in the symbolic dynamics of $f$). We have the following lemma:
 \begin{lemme}
  For all $p\in E$, $(\lambda_{p,n})$ is an increasing sequence bounded by $1$. Let $\lambda_p$ be its limit. Then: 
$$\sum_{p \in E} \lambda_p = 1.$$
   \end{lemme}
 \emph{Proof.} Write $f^n=(P_n Q_{1,n}+R_{1,n},P_n Q_{2,n}+R_{2,n})$. Recall that $I(f^n)$ is the intersection of $L_{\infty}$ with the zero set of $P_n$. By Proposition \ref{composition}:
  $$P_{n+1}=(P_n)^DP(Q_{1,n},Q_{2,n}).$$
Hence, $(\lambda_{p,n})$ is increasing since $(P_n)^D$ is a factor of $P_{n+1}$.
 
 Set $d_n=\deg (P_n)$ and $d'_n=\deg (Q_{i,n})$. We deduce from Proposition \ref{composition}:
  \begin{eqnarray*}
  d'_n = (d')^n \quad  \text{and} \quad   d_n  = D^n-(d')^n.
 \end{eqnarray*}
So, $\sum_{p\in E} \lambda_{p,n} =\frac{d_n}{D^n} \to 1$. This completes the proof.  $\Box$  \\
 
 \noindent \emph{Remarks} \begin{enumerate} \item In a way, the indeterminacy points of $f^n$ take asymptotically all the available degree, so they carry the main part of the dynamics near $L_{\infty}$ (cf. Proposition \ref{intersection}).
 \item The sequence $(\lambda_{p,n})_n$ can be strictly increasing as we will see in the last two examples of Section \ref{examples}. One can check that $(\lambda_{p,n})_n$ is strictly increasing after some rank if and only if $p$ is preperiodic.
  \end{enumerate}

Recall that, on $\mathbb{C}^2$, for $f$ in $\mathcal{G}$ of degree $D$, the sequence of positive functions $(u_n=\frac{1}{D^n} \log^+ \|f^n(z,w) \|)$ almost decreases (i.e. $(u_n+c_n)_n$ is decreasing for some sequence of constant $(c_n)_n$ decreasing to zero) to the \emph{Green function} $u$ of $f$ which is a potential of the Green current $T$. Furthermore, the function $\tilde{u}(z,w)=u(z,w)-\frac{1}{2}\log(|z|^2+|w|^2+1)$ is a bounded quasi-plurisubharmonic function on $\mathbb{C}^2$, thus it extends to $\mathbb{P}^2$, and this extension satisfies $dd^c\tilde{u}=T-\omega_{FS}$ where $\omega_{FS}$ is the Fubini-Study form on $\mathbb{P}^2$. 

We will see in Proposition \ref{intersection} that $\tilde{u}_{|L_\infty}$ is not identically equal to $-\infty$ so we can define the measure $m_{\infty}=T\wedge [L_\infty]$ which is the trace of the Green current at infinity. Since the sequence of functions $\tilde{u}_n(z,w):=u_n(z,w)-\frac{1}{2}\log(|z|^2+|w|^2+1)$ is almost decreasing, $m_{\infty}$ is the limit in the sense of current of the sequence $((dd^c\tilde{u}_n(z,w)+\omega_{FS})\wedge [L_\infty])$. In particular, we have $m_{\infty}=dd^c(\tilde{u}_{|L_\infty})+(\omega_{FS})_{|L_\infty}$. The next proposition shows that $m_{\infty}$ is a combination of Dirac masses at the points of $E$, with computable coefficients. 
 
 \begin{proposition}\label{intersection}
Let $f$ be in $\mathcal{G}$ and $\tilde{u}$ be as above. For $p\in E$, we denote by $[a_p:b_p:0]$ its homogeneous coordinates. Then:
\begin{eqnarray*}
\tilde{u}([z:w:0])=\log^+(\prod_{p \in E}|a_p w-b_p z|^{\lambda_p})-\frac{1}{2}\log(|z|^2+|w|^2).
\end{eqnarray*}
In particular, we have the formula:
\begin{eqnarray*}
 m_{\infty} =\sum_{p \in E} \lambda_p \delta_p
\end{eqnarray*} 
where $\delta_p$ is the Dirac mass at $p$.
 \end{proposition} 

\emph{Proof.} With the above notations, we have that in $\mathbb{C}^2$:
\begin{eqnarray*}
\tilde{u}_n(z,w) &=& \frac{1}{D^n}\log^+\|(P_nQ_{1,n}+R_{1,n})(z,w),(P_nQ_{2,n}+R_{2,n})(z,w)\| \\
                 & & -\frac{1}{2}\log(|z|^2+|w|^2+1). 
\end{eqnarray*}
So, first apart from the point of $E$, and hence everywhere on $L_\infty$ by semi-continuity, the extension is given by:
$$\tilde{u}_n([z:w:0])=\frac{1}{D^n} \log \| (P_n Q_{1,n})(z,w),(P_nQ_{2,n})(z,w)\|-\frac{1}{2}\log(|z|^2+|w|^2)$$
By definition of the $\lambda_{p,n}$ and Corollary \ref{stability}, there is a constant $C_n$ depending on the choice of the coordinates of the elements of $E$ such that:
 $$ \tilde{u}_n([z:w:0])=\sum_{p\in E} \lambda_{p,n} \log|a_p w-b_p z|+ \frac{1}{D^n} \log \|f_{\infty}^n[z:w]\|+C_n -\frac{1}{2}\log(|z|^2+|w|^2).$$
From one-dimensional theory, we know  that $\frac{1}{d'^n} \log |(f_{\infty})^n[z:w]-\frac{1}{2}\log(|z|^2+|w|^2)$ converges to a continuous function on $L_\infty$ and $\sum_{p\in E} \lambda_{p,n} \log(|a_p w-b_p z|)$ converges thanks to the previous lemma. The last identity and the fact that $d' < D$ imply the first formula in the proposition. The formula giving $m_{\infty}$ is then clear by the Poincaré formula. $\Box$ \\
  
 \noindent \emph{Remark.} The previous proof can be applied to all the algebraically stable polynomial maps of $\mathbb{C}^2$ with indeterminacy points on $L_{\infty}$.

  \subsection{Topological degree}\label{hypothese}
 Let $N$ be a small enough neighborhood of $L_{\infty}$ and $V$ be a neighboorhood of $I(f)$, then there are constants $C$ and $C'$ such that for $p$ in $N\backslash V$, we have:
\begin{equation*}
 C \|p\|^{D}\leq \|f(p)\|\leq C' \|p\|^{D}.
 \end{equation*}
 Let us assume here that the considered mapping satisfies in addition: for all $I \in I(f)$, there exist a number $l_I$, a neighborhood $V(I)$ of $I$, a neighborhood $V(f_{\infty}(I))$ of $f_{\infty}(I)$ and constants $C_1$ and $C_2$ such that for all $p \in V(I)$ with $f(p)\notin V(f_{\infty}(I))$, we have:
 \begin{equation}
 \label{condition}
 C_1 \|p\|^{l_I}\leq \|f(p)\|\leq C_2 \|p\|^{l_I}
 \end{equation}
This condition is easy to check in practice. Under these assumptions, we can compute the topological degree of $f$ which is the mass of the pull-back of any probability measure by $f$. The difference with the case with no dynamics on $L_{\infty}$ is that we have to count the number of preimages of a generic line by $f_{\infty}$. We have the following proposition:
 \begin{proposition}  
Let $f \in \mathcal{G}$ satisfying (\ref{condition}). Then the topological degree of $f$ is given by:
 $$d_t =\sum_{I\in I(f)} l_I \alpha_I+d'D .$$
In particular, we have $d_t>D$. 
  \end{proposition}
\emph{Proof.} Let $L$ be a generic line, we consider the probability measure $[L_\infty]\wedge[L]$ (which is the Dirac mass at the intersection of $L$ and $L_\infty$). By definition, its pull back by $f$ is of mass $d_t$. After some change of coordinates, we can assume that the point $[1:0:0]$ is not on $L$ and $f^{-1}(L)$. So we work in the coordinates $(u,v)=(Z/W,T/W)$ where a potential of $L_\infty=(v=0)$ is $\varphi(u,v)=\log|v|$. We must compute:
$$\int_{\mathbb{P}^2} f^{*} ([L_\infty]\wedge[L])= \int_{f^{-1}(L)} dd^c (\varphi \circ f).$$ 
For each $I$ in $I(f)$, let $\mathbb{B}_I$ be a bidisk in $V(I)$ for the $(u,v)$ coordinates, and for each $p$ in $f^{-1}_\infty(L_\infty \cap L)$ let  $\mathbb{B}_p$ be a bidisk around $p$. Since $L$ is a generic line, we can assume that $f^{-1}_\infty(L_\infty \cap L) \cap I(f)=\varnothing$  and that all those bidisks are disjoint. The previous integral become:
$$d_t= \sum_{I \in I(f)} \int_{f^{-1}(L)\cap \mathbb{B}_I} dd^c (\varphi \circ f) +\sum_{p \in f^{-1}_\infty(L_\infty \cap L)} \int_{f^{-1}(L)\cap \mathbb{B}_p} dd^c (\varphi \circ f).$$
Observe that $\varphi\circ f-l_I\log|v|$ is a bounded pluriharmonic function on $\mathbb{B}_I\backslash L_\infty$ thanks to (\ref{condition}), so it defines in fact a pluriharmonic function on $\mathbb{B}_I$. Hence, on these bidisks, $dd^c(\varphi\circ f)$ is equal to $l_I$ times the current of integration on $L_\infty$. Using the same argument for $\mathbb{B}_p$, we deduce:
$$d_t=\sum_{I \in I(f)} \int_{f^{-1}(L)\cap \mathbb{B}_I} l_I dd^c (\log|v|) +\sum_{p \in f^{-1}_\infty(L_\infty \cap L)} \int_{f^{-1}(L)\cap \mathbb{B}_p} D dd^c (\log|v|)$$
 which is what we wanted since $\int_{f^{-1}(L)\cap \mathbb{B}_I}  dd^c (\log|v|)=\alpha_I$ is the intersection multiplicity at $I$ of $L\infty$ and $f^{-1}(L)$ and since there are $d'$ preimages of $L\cap L_\infty$ by $f_\infty$. $\Box$

\section{Structure of the Julia set and of the Green current near infinity}
 We assume in this section that $f_{\infty}$ is uniformly hyperbolic (i.e. the forward orbit of each critical point converges towards some attracting periodic orbit), and that the indeterminacy points are not in the Julia set $J_{\infty}$ of $f_\infty$. After a unitary change of coordinates, we can assume that $[1:0:0]$ is not in $J_{\infty} \cup E$. Hence $(u,v)=(Z/W,T/W)$ is a coordinate system of a neighborhood of $L_\infty \backslash [1:0:0]$ where $L_\infty=(v=0)$. We also need the hypothesis that the indeterminacy points are not periodic.\\
 
We construct suitable boxes (polydisks) around the elements of $E$ such that $f$ defines horizontal-like maps between these boxes. 

 After decomposing the Julia set into pieces near infinity, we construct a subshift on $E^{\mathbb{N}}$. We then decompose the Green current along these pieces by pulling-back a smooth vertical positive closed $(1,1)$-form in the boxes which gives the Green current in a neighborhood of infinity. Observe that the set of maps we consider contains an open set in the space of parameters. 

Next, we give an application for the escape rate of $f$ and we explain how to obtain a weaker decomposition in the more general case where some indeterminacy points are in $J_\infty$. Finally we explain our results through examples. 

\subsection{Construction of the boxes}\label{construction}
The purpose of this section is to prove the following proposition:  
  \begin{proposition}\label{decoupage}
  For all $p$ in $E$, there is a bidisk $\Delta_p$ centered at $p$ such that $f$ induces by restriction a horizontal-like map from $\Delta_p$ to $\Delta_q$ for all $q\in E$ if $p \in I(f)$ and for $q=f_{\infty}(p)$ if $p$ is not an indeterminacy point. We denote by $f_{p,q}$ this restriction.
  
The bidisks can be taken arbitrarily small. We can choose them so that for all $I\in I(f)$ and all $q\in E -\{I\}$ then $\Delta_I \cap \Delta_q =\varnothing$,  and for all $p \in E$ and all $q,q' \in f_\infty^{-1}(p)$ then $\Delta_q\cap\Delta_{q'}=\varnothing$.
  \end{proposition}  
 
 Since $f_{\infty}$ is uniformly hyperbolic, we can put a smooth conformal metric $g$ on $L_{\infty}$ such that $\|Df_{\infty}(z)\|_{g} \geq \lambda>1$ on $J_{\infty}$. Let us remark that $E$ is discrete in the Fatou set $F_{\infty}:=L_\infty \backslash J_\infty$ since the only components of $F_{\infty}$ are basins of attraction of periodic cycles and $\overline{E}=E\cup J_{\infty}$ (see \cite{MIL}). The idea is first to construct disks on $L_{\infty}$ which will be thickened to get bidisks. So, we use the following lemma:
 \begin{lemme}\label{1dim}
There is a constant $c>0$ such that for all $p$ in $E$, there exists a disk $\mathbb{D}_p$ for the metric $g$ such that if $p\in I(f)$ and $q \in E$ then $\textrm{dist}_g(f_{\infty}(\partial \mathbb{D}_p), \mathbb{D}_q) \geq c$ and if $p\in E\backslash I(f)$ then $\textrm{dist}_g(f_{\infty}(\partial \mathbb{D}_p), \mathbb{D}_{f_{\infty}(p)}) \geq c$. Furthermore, we can choose the radii of those disks to be bounded and arbitrarily small. 
 \end{lemme} 
 \emph{Proof.} 
  Let $U$ be an open neighborhood of $J_{\infty}$ in $L_{\infty}$ with smooth boundary such that $\|Df_{\infty}(z)\|_g\geq \rho >1$ on $U$ and $f_{\infty}^{-1} U \subset U$. There is only a finite number of elements of $E$ in $L_{\infty} \backslash U$. Modifying $U$ if necessary, we can assume that $I(f) \cap U= \varnothing$ and $\partial U \cap E =\varnothing$.
  
  For $I \in I(f)$ such that $f_{\infty}(I) \notin E$, we consider $\mathbb{D}_I$ a disk centered at $I$ on $L_{\infty}$ for the metric $g$ with $f_{\infty}(\mathbb{D}_I)$ far from the other points of $E$. Restricting $\mathbb{D}_I$ if necessary, we can assume that for all $p$ in $f_{\infty}^{-1} \{ I \}$ there is a  disk $\mathbb{D}_p$ centered in $p$ on $L_{\infty}$ such that $f_{\infty}(\partial \mathbb{D}_p) \cap \mathbb{D}_I = \varnothing$ (we use the fact that $f_{\infty}$ is open). We iterate this construction with the preimages of all the $p$ till all of them are in $U$. Of course, we may have to shrink $\mathbb{D}_I$ at each step. We apply this process to all the elements of $I(f)$ such that $f_{\infty}(I) \notin E$.
   
    Since we assumed the elements of $I(f)$ are not periodic for $f_{\infty}$, we then have disks $\mathbb{D}_p$ for all the $p$ in $ E \backslash U $ such that $f_{\infty}(\partial \mathbb{D}_p) \cap \mathbb{D}_{f_{\infty}(p)} = \varnothing$. Let $r$ be the smallest radius of all these disks. It can be chosen arbitrarily small.
    
   Next, by hyperbolicity, there is some $\varepsilon_0 >0$ such that $f_{\infty}$ is injective on any disk $\mathbb{D}_g(z,\varepsilon_0)$ for all $z$ in $f_{\infty}^{-1} U$, and is closed to its differential. Namely, for all $\varepsilon \leq \varepsilon_0$, there is a $\rho'>1$ such that we have $\mathbb{D}_g(f(z),\rho' \varepsilon) \Subset f_\infty(\mathbb{D}_g(z,\varepsilon))$ for $z$ in $f_{\infty}^{-1} U$. Then, for $r$ small enough, we have some $r'$ such that for all $p$ in $E \cap f_{\infty}^{-1} U$, the disk $\mathbb{D}_p=\mathbb{D}_g(p,r')$ satisfies $f_{\infty}(\partial \mathbb{D}_p) \cap \mathbb{D}_{f_{\infty}(p)} = \varnothing$. The existence of the constant $c>0$ is then clear by construction for $p\in L_{\infty} \backslash f^{-1}_\infty(U)$ and by hyperbolicity for $p\in E\cap f^{-1}_\infty(U)$. $\Box$\\
   
\noindent\emph{Proof of Proposition \ref{decoupage}.}
Recall that the line at infinity is $f$-attracting: there is a constant $C>1$ such that for $M=(z,w)$ in $\mathbb{C}^2$ with $\|M\| \geq A$, we have $\|f(M)\| \geq C \|M\|$ where $\| \star \|$ denotes the euclidean norm. Furthermore:
 $$\| M \|^2=|z|^2+|w|^2=|\frac{u}{v}|^2+|\frac{1}{v}|^2.$$
If $p=(u_p,v_p)$, define $\Delta_p=\mathbb{D}_p \times \{ v< \frac{\epsilon}{\sqrt{1+|u_p|^2}} \}$ with $\varepsilon$ small. For $M=(u,v) \in \Delta_p$, we have that $$(1+\nu)^{-1}(|u_p|^2+1)|\frac{1}{v}|^2 \leq \| M \|^2\leq (1+\nu)(|u_p|^2+1)|\frac{1}{v}|^2$$
where $\nu>0$ depends only on the radius of $\mathbb{D}_p$ (since $u$ is uniformly bounded) and goes to zero with it.  We take the radii of the $\mathbb{D}_p$ small enough so that the square of the norm of an element in $\Delta_p$ is close to $(|u_p|^2+1)|\frac{1}{v}|^2$.\\
 
We choose $\varepsilon$ so that all the $\Delta_p$ are in the domain where the infinity is attracting. Restricting $r'$ which is the supremum of the radii of the disks $\mathbb{D}_p$ if necessary, we can assume that $f(\Delta_p) \cap \partial \Delta_q \subset \partial_v \Delta_q$ for all $q$ if $p$ is an indeterminacy point and for $q=f(p)$ otherwise. Now, using the uniform continuity of $f$ in the complement of some small neighborhood of the indeterminacy set and the existence of $c$ in Lemma \ref{1dim}, we can choose $\varepsilon$ so that $f(\partial_v \Delta_p) \cap \Delta_{q}=\varnothing$ for all $q\in E$ if $p\in I(f)$ and for $q=f(p)$ otherwise. Finally, since the image of any small neighborhood of an indeterminacy point by $f$ contains the whole line at infinity, we have $f(\Delta_p)\cap\Delta_q \neq \varnothing$. The last part of the proposition is clear for $r'$ small enough (we use the hyperbolicity of $f$ once again here). $\Box$

\subsection{Construction of the subshift}\label{subshift}
Now, we define the symbolic dynamics which will appear in the decomposition of the Green current. First, we will need to know the degree of the horizontal-like maps $(f_{p,q})$. Recall that $\alpha_i$ is the multiplicity of the indeterminacy point $I_i \in I(f)$. We take the notations of $(\ref{form})$. The following lemma is clear:
 \begin{lemme}   
 \begin{enumerate}
   \item If $p$ is in $E \backslash I(f)$, then the degree of $f_{p,q}$ is the local degree of $f_{\infty}$ at $p$,  
   \item if $p=I_j$ is in $I(f)$ and $q \neq f_{\infty}(I_j)$, then the degree of $f_{p,q}$ is $\alpha_j$,
   \item if $p=I_j$ is in $I(f)$ and $q = f_{\infty}(I_j)$, then the degree of $f_{I_j,f_{\infty}(I_j)}$ is the sum of $\alpha_j$ and the local degree of $f_{\infty}$ at $I_j$.
  \end{enumerate}
 \end{lemme}
Define $\Sigma' = E^{\mathbb{N}}$ and $\Sigma = \{ (\alpha_n)\in \Sigma', f_{u_n,u_{n+1}} \ \mathrm{exists} \}$,
the space of itineraries between indeterminacy points and theirs preimages. 
We consider the left shift $\sigma$ on $\Sigma$ and $\Sigma'$. Define $N=\bigcup_{p\in E} \Delta_p $. For $\alpha \in \Sigma$, let us consider:
 $$\mathcal{K}_{\alpha}=\{ p \in  N, \ f^j(p) \in \Delta_{\alpha(j)} \}. $$
 Then, for all $\alpha \in \Sigma$, $\overline{\mathcal{K}_{\alpha}}$ is not empty as a decreasing intersection of vertical closed sets in $\Delta_{\alpha(0)}$. Let  $\mathcal{K}$ be the union of all the $\mathcal{K}_{\alpha}$ so that $\mathcal{K} \subset N$.\\
 
Observe that $T\wedge[L_\infty]$ is the slice of $T$ by $(v=0)$. Using the formula giving the trace of $T$ on $L_{\infty}$ and the invariance of $T$ ($f^*T=DT$), we have that:
 $$ \forall p \in E, \lambda_p= \frac{1}{D} \sum_{q\in E} d_{p,q} \lambda_q $$ 
with the convention that $d_{p,q}=0$ if $f_{p,q}$ is not defined. For all $p\in E$, we deduce:
\begin{eqnarray}\label{sum}
1 = \sum_{q\in E} \frac{d_{p,q}\lambda_q}{D\lambda_p}.
\end{eqnarray}
For example, if $p$ is not an indeterminacy point, we have that all the $d_{p,q}$ are zero except for $q=f(p)$ and the formula becomes:
   $$1 = \frac{d_{p,f(p)}\lambda_{f(p)}}{D\lambda_p}.$$
And if $p=I$ is in the indeterminacy set with $d_{I,q}$ constant (i.e. $f_{\infty}(I) \notin E$), then:
   $$ 1 = \sum_{q\in E} \lambda_q.$$  
Let $A:=(a^q_p)_{p,q \in E}$ be the infinite matrix defined by $a^q_p= \frac{d_{p,q}\lambda_q}{D\lambda_p}$. The entry $a^p_q$ can be seen as the \emph{probability to go from $\Delta_p$ to $\Delta_q$ by $f$} in term of slice mass (see the proof of Theorem \ref{main}). Of course, if $p$ is not an indeterminacy point, one always goes to $\Delta_{f(p)}$ (the probability is 1). We put on $\Sigma$ the Borel measure $\nu$ defined by:
 $$\nu( \{ \alpha \in \Sigma, \ \alpha(0)=\alpha_0, \dots , \alpha(n)=\alpha_n \})=\lambda_{\alpha_0} \times \prod^{n-1}_{i=0} a^{\alpha_i+1}_{\alpha_{i}}= \lambda_{\alpha_n} \times \prod^{n-1}_{i=0} \frac{d_{\alpha_i,\alpha_{i+1}}}{D}.$$  
\begin{proposition}
 The left shift $\sigma$ on $\Sigma$ defines a subshift for which the measure $\nu$ is invariant and mixing. 
\end{proposition}
\emph{Proof.} Definitions and facts on symbolic dynamics and especially subshift can be found in \cite{K.H}, pp. 156-158. By (\ref{sum}), we already have for all $p\in E$ that 
\begin{eqnarray}
\sum_q a_p^q=1
\end{eqnarray}
What remains to be proved is that the vector $(\lambda_p)$ is an eigenvector for the matrix $^tA$ associated with the eigenvalue 1 (that gives the invariance of $\nu$). That is: 
\begin{eqnarray}
\sum_{p} a^q_p \lambda_p= \lambda_q
\end{eqnarray}
 which is clear again by (\ref{sum}). \\

Furthermore, the matrix $A$ is transitive in the sense that for each $(p,q)$ the entry of index $(p,q)$ in $A^n$ is stricly positive for some $n$ (it is clear if $p$ is in the indeterminacy set for $n=1$ and if $p \in f^{-j}(I(f))$, then it is true for $n=j+1$). We deduce that $\nu$ is mixing. $\Box$\\

 \noindent \emph{Remark.} We can consider only a finite part of $E$ containing the indeterminacy points and their preimages up to some order and regroup the rest of the elements of $E$ in a same box. Then we can obtain a finite Markov chain, but we lose some part of the information. 
 
\subsection{Decomposition of the Green current}
 Let us denote by $\mathcal{L}_{p,q}$ the operator $\frac{1}{d_{p,q}}f^{*}_{p,q}$ acting on vertical currents. Recall that $f$ is a polynomial map of $\mathbb{C}^2$ having indeterminacy points on $L_\infty$ which is $f$-attracting. The map $f_\infty$ is hyperbolic and the indeterminacy points of $f$ are on the Fatou set of $f_\infty$ and not periodic. We can now prove our main theorem:
\begin{theorem}\label{main}

1. There exists an at most countable set $\Theta \subset \Sigma$ such that for all $\alpha\in \Sigma \backslash \Theta$, there is a unique current $T_\alpha$ satisfying the following property: for all sequence of currents $(S_{k+1})$ of bidegree (1,1), positive, closed, vertical in $\Delta_{\alpha(k+1)}$ of slice mass $1$, we have:
 $$\mathcal{L}_{\alpha(0),\alpha(1)} \dots  \mathcal{L}_{\alpha(k),\alpha(k+1)}S_{k+1}\to T_\alpha.$$ 
2. The Green current $T$ admits the following decomposition in $N$:
 $$T=\int_{\Sigma} T_\alpha d\nu(\alpha).$$
 \end{theorem} 
\emph{Proof.} Since $\mathcal{K}=\bigcup \mathcal{K}_{\alpha}$, only a countable number of $\mathcal{K}_{\alpha}$ have positive Lebesgues measure. Then Theorem \ref{convergence} implies the first part. 
 
 For the second part, let $S'_p$ be a smooth positive closed (1,1)-form in $\mathbb{P}^2$ such that near $L_\infty$, $S'_p$ has its support in $\Delta'_p=\mathbb{D}'_p \times \{ v< \frac{\epsilon}{\sqrt{1+|u_p|^2}} \}$, with $\mathbb{D}'_p \Subset \mathbb{D}_p$. Let $S_p$ be the restriction of $S'_p$ to $\Delta_p$, it is a vertical positive  closed current. Normalize $S'_p$ so that $S_p$ is of slice mass $\lambda_p$. Observe that if $S'=\sum S'_p$ then $\lim \frac{1}{D^n}(f^n)^* (S')=\sum \lim  \frac{1}{D^n}(f^n)^* (S'_p)=T$ since $\sum \lambda_p=1$. Define $S=\sum S_p$.  Finally, write:
  $$\Sigma_n=\{(a_0,a_1,\dots,a_{n-1})\in E^n | \exists \alpha \in \Sigma, \forall i \leq n-1, \ a_i=\alpha(i)\}$$
 and for $a\in \Sigma_n$, write $C_a$ for the cylinder: 
 $$\{ \alpha \in \Sigma| \ \alpha(0)=a_0,\dots ,\alpha(n-1)=a_{n-1} \}.$$
Pulling back $S$ by $f$ gives:
\begin{eqnarray*}
\frac{1}{D} f^{*} S & = & \frac{1}{D}\sum_{p,q\in E}  f^{*}_{p,q}S_q \\
                    & = & \sum_{\alpha \in \Sigma_1 } \frac{d_{\alpha(0),\alpha(1)}}{D} \mathcal{L}_{\alpha(0),\alpha(1)}S_{\alpha(1)}.
\end{eqnarray*}
We iterate:
 \begin{eqnarray*}
\frac{1}{D^k} (f^k)^* S & = & \sum_{\alpha\in \Sigma_k}  \frac{\prod^{k-1}_{i=0} d_{\alpha(i),\alpha(i+1)}}{D^k} \mathcal{L}_{\alpha(0),\alpha(1)}\dots\mathcal{L}_{\alpha(k-1),\alpha(k)} S_{\alpha(k)} \\
                        & = & \sum_{\alpha \in \Sigma_k } \nu(C_{\alpha}) \mathcal{L}_{\alpha(0),\alpha(1)}\dots\mathcal{L}_{\alpha(k-1),\alpha(k)}  \frac{S_{\alpha(k)}}{\lambda_k}.
 \end{eqnarray*}
The left hand side goes to the Green current $T$. By the first part of the theorem, the general term of the right hand side tends to $T_\alpha$ for $\alpha$ generic so we get the result by dominated convergence. $\Box$\\
 
 \noindent \emph{Remark.} The dynamics of $f$ near infinity is semi-conjugated to the subshift $\sigma$ in the sense that $f(\mathcal{K}_\alpha) \subset \mathcal{K}_{\sigma(\alpha)}$.\\  
 We see that the current is carried by $\mathcal{K}$ which does not meet $J_\infty$. So, as announced in the introduction, the local stable manifolds to the Julia set of $f_{\infty}$ do not carry any part of the Green current, but they are contained in its support.

\subsection{Escape rate} 
We take $f\in\mathcal{G}$ satisfying the condition (\ref{condition}), we also suppose that $f_{\infty}(I(f)) \cap E =\varnothing$ (else, we would have $p$ in $E$ such that $f(p) \in f(V(I(f)))$). We want to compute the possible values of the \emph{upper escape rate} $\bar{l}$ where $\log(\bar{l})=\lim\sup \frac{1}{n}  \log^+ \log^+ \| f^n \|$ which becomes $\lim\sup \frac{1}{n}  \log \log (\| f^n \|)$ in $N$. In the same way, we define the \emph{lower escape rate} $\underline{l}$ and we are interested in knowing where these two functions match up, in which case we note $l$ their common value which we simply call the \emph{escape rate}.\\
 
 For $p\in E\backslash I(f)$, we set $l_p=D$. We have the following lemma:
 \begin{lemme} 
Let $\alpha\in\Sigma$ and $q\in \mathcal{K}_\alpha$. We have:
$$\frac{1}{n} \log \log \| f^n(q) \| = \frac{1}{n} \log(l_{\alpha(0)}l_{\alpha(1)}\dots l_{\alpha(n-1)})+O(\frac{\log n}{n}).$$ 
 \end{lemme}
 \emph{Proof} We have constants $c_1$ and $c_2$ such that: 
$$ c_1 \leq \log \|f^{j+1}(q)\| - l_{\alpha(j)}\log \|f^j(q)\| \leq c_2.$$
Taking a combination of these inequalities for $j\leq n-1$ gives:
 
  \begin{eqnarray*}
c_1 \Big( \sum_{j=0}^{n-1} l_{\alpha(j+1)}\dots l_{\alpha(n-1)}\Big) +l_{\alpha(0)}\dots l_{\alpha(n-1)}\log\|q\| \leq \log\|f^n(q)\|, 
 \end{eqnarray*}
with a similar inequality for the right hand side. Taking the logarithm and dividing by $n$ give:
 $$\Big| \frac{1}{n} \log \log \| f^n(q) \|- \frac{1}{n} \log(l_{\alpha(0)}\dots l_{\alpha(n-1)}) \Big| \leq \frac{1}{n}\log \Big(\log \|q\|+C \sum_{j=0}^{n-1} \frac{1}{l_{\alpha(0)}\dots l_{\alpha(j-1)}} \Big) $$
The sum in the right hand side is a $O(n)$ which concludes the proof. $\Box$\\
 
 Choosing a suitable $\alpha$, we deduce from the lemma that the range of the escape rate in $N$ is $[ \min l_I, D]$ (the details are left to the reader). In this case, it is interesting to observe that the set of possible escape rates is an interval which is a new property for polynomial mappings. Let $\lambda$ denote  the slice mass $1-\sum_{I\in I(f)} \lambda_I$ of $T$ outside a neighborhood of $I(f)$. We have the following theorem:
 \begin{theorem}
 For $\|T\|$-almost every point $q$ in $N$, the escape rate $l(q)$ exists and is equal to $D^{\lambda}\prod_{I\in I(f)} l_I^{\lambda_I}$.
  \end{theorem}
 \emph{Proof.}
  Since the left shift $\sigma$ is ergodic for $\nu$, the Birkhoff's ergodic theorem yields that for $\nu$-almost every $\alpha$:
  $$\exp \Big( \frac{1}{n}\sum_{i=0}^{n-1} \log l_{\sigma^i(\alpha)(0)} \Big) \to \exp(\int_{\Sigma} \log l_{\alpha(0)} d\nu)  = D^{\lambda}\prod_{I\in I(f)} l_I^{\lambda_I}.$$
And the theorem follows from the previous lemma and Theorem \ref{main}. $\Box$ 
 
\subsection{Generalization}\label{generalization} 

In the case where some indeterminacy points are on $J_\infty$ (possibly periodic), we can obtain a decomposition of the Green current by building a cover of $J_\infty$ by disks such that for all $\mathbb{D}$ in this cover, there exist disjoint disks $\mathbb{D}_1, \mathbb{D}_2, \dots, \mathbb{D}_{d'}$ in the cover with $f_\infty^{-1}(\mathbb{D}) \subset \mathbb{D}_1 \cup \mathbb{D}_2\cup \dots \cup \mathbb{D}_{d'}$ and $\mathbb{D} \Subset f_\infty (\mathbb{D}_i)$ for all $i\leq d'$. The trick is to have two disks $\mathbb{D}_I$ and $\mathbb{D}'_I$ around each indeterminacy point $I\in I(f)$ so that $\partial f_\infty(\mathbb{D}_I) \cap \overline{\mathbb{D}}=\varnothing$ or $\partial f_\infty(\mathbb{D}'_I) \cap \overline{\mathbb{D}}=\varnothing$. Finally, we follow the construction of Section \ref{construction} with $U$ being replaced by the union of all those disks.

This time we only have a finite number of bidisks and when we pull back the Green current near some point of $E$ to an indeterminacy point $I$ in $J_\infty$, we may have to choose between the two bidisks centered at $I$ in order to have a horizontal-like map. We only get a finite subshift, but taking a finer cover, we get more precision on the decomposition (only on a smaller neighborhood of $L_\infty$). Somehow the decomposition is not intrinsic because we do not pull back according to the itinerary but it assures that the Green current is not extremal in a neighborhood of $L_\infty$.

\subsection{Examples}\label{examples}
First let us explain our results in two examples where the dynamics at infinity is linear.\\

\noindent \emph{Example 1.} Consider the case where $f_\infty$ is given by $u\mapsto 2u$ and where the indeterminacy set is reduced to $(1,0)$ with multiplicity 1 in the $(u,v)$ coordinates (thanks to Proposition \ref{criterion}, we know this case exists, take for example $f(z,w)=C(2z(z-w)+z,w(z-w))$ for $C$ large enough). Then, using Proposition \ref{intersection}, we find that:
\begin{itemize}
\item $E=\{ p_n=(\frac{1}{2^n},0) \ n\geq0 \}.$ 
\item $\lambda_{n} = \frac{1}{2^n}$
\item the matrix of the subshift is:
$$
\begin{pmatrix}
         \frac{1}{2} & \frac{1}{4}  & \frac{1}{8}    & \dots   \\
         1           & 0            & 0              & \dots   \\              
         0           & 1            & 0              & \dots   \\  
         \vdots      & \vdots       & \vdots         & \ddots

\end{pmatrix}
$$
\end{itemize}
  An element $\alpha\in\Sigma$ can be written $(p_{n_1}, p_{n_1-1},\dots,p_0,p_{n_2},\dots,p_0,\dots)$ for some sequence $(n_i)$ in $\mathbb{N}$. The dynamics in the space of itineraries is simple: a point in $\mathcal{K}_\alpha$ where $\alpha_0=p_{n_1}$ is sent near $p_{n_1-1}$ then near $p_{n_1-2}$ and so on untill it arrives near $p_0$, in which case it can be sent near any element of $E$ since $p_0$ is an indeterminacy point.\\

\noindent \emph{Example 2.} This time, we still take $f_\infty$ given by $u\mapsto 2u$ and we suppose that the indeterminacy points are $I_0=(2,0)$ and $I_1=(1,0)$ with multiplicity 1 in the $(u,v)$ coordinates, so $D=3$ (for example: $f(z,w)=(2z(z-w)(z-2w)+z^2,w(z-w)(z-2w))$). In this case, we have that $f_\infty^{-1}I_0=I_1$. Again, using Proposition \ref{intersection}, we find that:  
  \begin{itemize}
\item $E=\{ p_n=(\frac{1}{2^{n-1}},0) \ n\geq0 \}$.
\item We have $\lambda_0=\lambda_{I_0}=\frac{1}{3}$, $\lambda_1=\lambda_{I_1}=\frac{4}{9}$, $\lambda_{p_n}=\frac{4}{2^{n+1}}$. 
\item The matrix of the subshift is:
$$
\begin{pmatrix}
         \frac{1}{3} & \frac{4}{9}  & \frac{4}{27}    &  \frac{4}{3^4}      &\dots       \\
         \frac{1}{2} & \frac{1}{3}  & \frac{1}{9}     &  \frac{1}{3^3}      &\dots       \\              
         0           & 1            & 0               &  0                  &\dots       \\  
         0           & 0            & 1               &  0                  &\dots       \\
         \vdots      & \vdots       & \vdots          &  \vdots             &\ddots
                 
\end{pmatrix}
$$
\end{itemize}

 The interesting fact here is that the entries of the second row are not proportionnal to the the slice mass, indeed a point near $I_1$ will have "more chances" to be sent on $\Delta_{I_0}$ by $f$ since $f_\infty(I_1)=I_0$.\\
 
Now, we consider the case were the indeterminacy points are in the exceptionnal set of $f_{\infty}$ (namely $f^{-1}(I)=I$). Observe that this case does not satysfy the hypothesis of Theorem \ref{main} since the indeterminacy points are periodic.\\

\noindent \emph{Example 3.}  The method used to treat the example easily extends to the case where $f_{\infty}$ admits a totally invariant point (i.e. $f_{\infty}^{-1}(e)=e$) which is equal to the indeterminacy set but for the sake of simplicity,we only consider:
$$f : (z,w) \mapsto (z^3+w^2, zw^2).$$
By Proposition \ref{criterion}, $L_{\infty}$ is $f$-attracting. We even have $\|f(z,w) \| \geq \|(z,w)\| ^2$ for $\|(z,w)\|$ large enough. The meromorphic extension of $f$ to $\mathbb{P}^2$ is given by: $f([Z:W:T])=[Z^3+ TW^2:ZW^2:T^3]$. The indeterminacy set of $f$ is reduced to $I_0=[0:1:0]$ and the dynamics at infinity is given by $f_\infty:[z:w:0] \mapsto [z^2:w^2:0]$ (so $f^{-1}_\infty(I_0)=I_0$). Thus $f$ is in $\mathcal{G}$ and is algebraically stable.
  
   The topological degree $d_t$ of $f$, which is by definition the number of preimages of a generic point, is equal to 8 (solve $f(z,w)=(0,1)$). It is greater than the algebraic degree.
  
   We use the coordinates $(u,v)=(\frac{Z}{W},\frac{T}{W})$ in which $L_\infty$ is given by $(v=0)$. The map $f$ becomes:
 $$f:(u,v) \mapsto \big( \frac{u^3+v}{u},\frac{v^3}{u} \big).$$ 
In these coordinates, the point $I_0$ becomes $(0,0)$. The map $f_\infty$ is given by $u \mapsto u^2$ for which the Julia set $J_\infty$ is the unit circle $(|u|=1)$. We have the following lemma: 
 \begin{lemme}
Let $V=\{(u,v),\ |u| < \frac{1}{2} \ \text{and} \ |v| < \frac{1}{4} |u|^3 \}$, then $f(V)\subset V$.
 \end{lemme}
\emph{Proof.} Let $(u,v)$ be in $V$. We check:
 $$\frac{|u^3+v|}{|u|}\leq |u|^2+\frac{|v|}{|u|} < \frac{1}{4}+\frac{|u|^2}{4} < \frac{1}{2}$$
We also have the inequalities:
 \begin{eqnarray*}
 \frac{|u^3+v|}{|u|} & \geq & |u|^2 - \frac{|v|}{|u|} >  |u|^2-\frac{|u|^2}{4} > \frac{1}{2}|u|^2 \\
 \frac{|v|^3}{|u|}   & < & \frac{1}{4^3} |u|^8.
 \end{eqnarray*}
 It is then sufficient to check that:
  $$ \frac{1}{4^3} |u|^8 < \frac{1}{4} (\frac{1}{2}|u|^2)^3$$
  which is obvious. $\Box$ \\
  
 We deduce from the lemma that $V$ is in the Fatou set since the sequence of iterates is normal there. Let then $\mathbb{D}_0 \subset \mathbb{D}_1$ be disks on $L_{\infty}$ centered on $I_0$, small enough to be contained in $V$, with $f_\infty^{-1}(\mathbb{D}_0) \Subset \mathbb{D}_1 $. Let $\mathbb{D}_2$ be a disk centered on $[1:0:0]$ containing the Julia set of $f_\infty$ with $\partial \mathbb{D}_2 \subset V$. We have that $f^{-1}(\mathbb{D}_2) \Subset \mathbb{D}_2$. We can shrink those disks to have $\mathbb{D}_1 \cap \mathbb{D}_2 = \varnothing$. \\
 
As in Proposition \ref{decoupage}, we want to "thicken" those disks in order to have bidisks such that $f$ defines by restriction horizontal-like maps between them. Close to $I$, the norm of a point (in the $(z,w)$ coordinates) is given by $|v|^{-1}$, but next to $[1:0:0]$, it is controled by $\frac{|u|}{|v|}$ so we use the coordinates $(u',v')=(\frac{T}{Z},\frac{W}{Z})$ there. Then, we define $\Delta_0 = \mathbb{D}_0 \times (|v|<\varepsilon)$, $\Delta_1 = \mathbb{D}_1 \times (|v|<\varepsilon)$ and $\Delta_2 = \mathbb{D}_2 \times (|u'|<\varepsilon')$. Take $\varepsilon$ and $\varepsilon'$ small enough so that the vertical boundaries of the bidisks are relatively compact in $V$. Observe that $\Delta_1 \backslash \Delta_0 \subset V$ is in the Fatou set of $f$.\\
 
 Recall that since $I_0$ is an indeterminacy point, any neighborhood of $I_0$ is sent on the whole $L_{\infty}$. Since $L_{\infty}$ is $f$-attracting, and by uniform continuity of $f$ away from any neighborhood of $I_0$, we can chose $\varepsilon$ and $\varepsilon'$ small enough so that:
 \begin{itemize}
 \item $f:\Delta_1 \to \Delta_0$ defines a horizontal-like map of degree 3 denoted by $f_{1,0}$
 \item $f:\Delta_1 \to \Delta_2$ defines a horizontal-like map of degree 1 denoted by $f_{1,2}$
 \item $f:\Delta_2 \to \Delta_2$ defines a horizontal-like map of degree 2 denoted by $f_{2,2}$ 
 \end{itemize}
Next, we consider the Green current $T$ of $f$. We know that its support is contained in the Julia set of $f$ (see \cite{SIB}). So we know that in some neighborhood of infinity, $T$ can be written as $T_1+T_2$ where $T_1$ and $T_2$ are vertical positive closed currents in $\Delta_0 \subset\Delta_1$ and in $\Delta_2$. Pulling-back $T_1$ and $T_2$ and using the invariance of $T$, we see that:
 $$ \frac{1}{3} f^*T=T=T_1+T_2$$
 So:
 \begin{eqnarray*}
 T_1 & = & \frac{1}{3} f^*_{1,0}T_1+\frac{1}{3}f^*_{1,2}T_2 \\
 T_2 & = & \frac{1}{3} f^*_{2,2}T_2
 \end{eqnarray*}
 Calling $m_1$ and $m_2$ the slice masses of $T_1$ and $T_2$, we can compute them using the previous equation and the fact that the pull-back of a vertical current of slice mass $m$ by a horizontal-like map of degree $d$ is of slice mass $dm$. So, we have:
 \begin{eqnarray*}
 m_1 & = & m_1 + \frac{1}{3} m_2 \\
 m_2 & = & \frac{2}{3} m_2
 \end{eqnarray*}
Hence, $m_2=0$ and so $T_2=0$. In particular, the support of the Green current of $f$ is \emph{strictly contained} in the Julia set $J$ since the stable manifolds associated to the Julia set $J_{\infty}$ of $f_{\infty}$ are in $J$ but $\text{supp}(T)$ does not meet $J_{\infty}$. In \cite{F.S}, there is a different example of such phenomenon.\\
For $\varepsilon>0$, we consider the small perturbation $f_\varepsilon$ defined by:
 $$f_\varepsilon : (z,w) \mapsto ((z+\varepsilon w)z^2+w^2, (z+\varepsilon w) w^2).$$
We check that $f_\varepsilon$ gives the same map at infinity than $f$ and that the indeterminacy point is now $I_\varepsilon=[-\varepsilon:1:0]$. We see that the preimages of $I_\varepsilon$ accumulate on the Julia set of $f_{\varepsilon \infty}$ and we have seen that they are on the support of the Green current which contrasts with what happens for $f$.\\

\noindent \emph{Example 4.} Consider the following families of polynomial maps:
\begin{eqnarray*}
 f(z,w) & = & \big(z^{n_1}w^{n_2}z^n+R_1(z,w), z^{n_1}w^{n_2}w^n+R_2(z,w)\big) \\
 g(z,w) & = & \big(z^{n'_1}w^{n'_2}w^{n'}+R'_1(z,w), z^{n'_1}w^{n'_2}z^{n'}+R'_2(z,w)\big)
\end{eqnarray*}
where the $n_i$, $n'_i$, $n$ and $n'$
are  positive integers and the $R_i$  and $R'_i$ are polynomials of degree smaller than $n_1+n_2+n$ and $n'_1+n'_2+n'$ respectively chosen so that $L_{\infty}$ is attracting for these polynomial mappings (use Proposition \ref{criterion}). The indeterminacy set is then $I=I(f)=I(g)=\{ [0:1:0], [1:0:0] \}$ which is totally invariant by $f_{\infty}$ or $g_{\infty}$ ($I$ is here the exceptional set of $f_\infty$ and $g_\infty$). Let us focus on the first case, as both cases can be understood through the same method.\\

As in the previous example, we can decompose the Green current in a neighboorhood of $L_{\infty}$ into $T_0$ and $T_1$ both being vertical in some bidisks $\Delta_0$ near $I_0=[0:1:0]$ and $\Delta_1$ near $I_1=[1:0:0]$. For $\alpha \in \{0,1\} ^{\mathbb{N}}$, define $\mathcal{K}_{\alpha} := \{ p \in \mathbb{C}^2, \ \text{s.t} \ \forall k \in \mathbb{N}, \ f^k(p) \in \Delta_{\alpha_k} \}\cup \{I_{\alpha(0)}\}$, which is non-empty since the image of any neighborhood of $I_i$ contains $L_{\infty}$. Let $\mathcal{K}$ be the union of all those sets. Define the matrix $A=(a_{ij})$ by:

\begin{displaymath}
A=\left( \begin{array}{cc}
          \frac{n_1+n}{n_1+n_2+n} &  \frac{n_2}{n_1+n_2+n} \\
          \frac{n_1}{n_1+n_2+n}   &  \frac{n_2+n}{n_1+n_2+n}
       \end{array} \right)
\end{displaymath}
We let $\lambda_0=\frac{n_1}{n_1+n_2}$ and $\lambda_1=\frac{n_2}{n_1+n_2}$. Then, we define the Borel measure $\nu$ on $\{ 0,1 \}^{\mathbb{N}}$ by:
 $$\nu(\{ \alpha \in \{ 0,1 \}^{\mathbb{N}}, \ \alpha(0)=\alpha_0, \dots , \alpha(n)=\alpha_n \})=\lambda_{\alpha_0} \times \prod^{n-1}_{i=0} a_{\alpha_{i}\alpha_{i+1}}.$$
We know that $\nu$ is invariant and mixing for the left shift $\sigma$ which defines a subshift on $\{ 0,1 \} ^{\mathbb{N}}$ (see Section \ref{subshift}). As in Theorem \ref{main}, we pull back $T_0$ and $T_1$ using the horizontal-like maps defined by restricting $f$ to the polydisks (in fact, we must take a large $\Delta_i'$ containing $\Delta_i$ and the unit circle $|u|=1$ and use the method of Section \ref{generalization}). Iterating the process, we obtain the following decomposition:\\
 
\noindent \emph{There exists an at most countable set $\Theta \subset \Sigma$ such that for all $\alpha$ in $\Sigma \backslash \Theta$, there exists a current $T_\alpha$ supported on $\overline{\mathcal{K}}_{\alpha}$ of slice mass 1 such that the Green current $T$ of $f$ admits the following decomposition:
 $$ T= \int_{\Sigma} T_{\alpha} d\nu(\alpha)$$} 
\bibliography{bible} 
 \noindent Gabriel Vigny, Mathématiques - Bât. 425, UMR 8628,\\ 
 Université Paris-Sud, 91405 Orsay, France.  \\
 \noindent Email: gabriel.vigny@math.u-psud.fr

\end{document}